\documentclass[a4paper]{article}  	
\usepackage{geometry}                		
\usepackage{graphicx}					
\usepackage{amssymb}
\usepackage{indentfirst}
\usepackage{amsmath}
\usepackage{amsthm}
\usepackage{subfigure}
\usepackage{siunitx}
\usepackage[square,sort,comma]{natbib}

\newtheorem{theorem}{Theorem}
\newtheorem{lemma}{Lemma}
\newcommand{\E}{\mathrm{E}}
\newcommand{\Var}{\mathrm{Var}}
\newcommand{\Cov}{\mathrm{Cov}}

\newcommand{\tr}{\mathrm{tr}}
\newcommand{\BigO}[1]{{\rm O}\left(#1\right)}

\newcommand{\iid}{\textrm{i.i.d.\ }}

\providecommand{\keywords}[1]{\textbf{\textit{Keywords: }} #1}

\title{V-Spline and Bayes Estimate}

\author{Zhanglong Cao, David Bryant, Matthew Parry}

\date{}  
\begin{document}
\maketitle

\begin{abstract}
It is known that a smoothing spline can be thought of as the posterior mean of a Gaussian process regression in a certain limit. By constructing a reproducing kernel Hilbert space with an appropriate inner product, the Bayesian form of the V-spline is derived when the penalty term is a fixed constant instead of a function. An extension to the usual generalized cross-validation formula is utilized to find the optimal V-spline parameters. 
\end{abstract}

\keywords{V-spline, generalised cross-validation, Bayes estimate, reproducing kernel Hilbert space}

\section{Introduction}

A Hilbert space is a real or complex inner product space with respect to the distance function induced by the inner product \citep{dieudonne2013foundations}. In particular, the Hilbert space $\mathcal{L}_2[0,1]$ is the set of square integrable functions $f(t):[0,1]\mapsto \mathbb{R}$, where all functions satisfy 
\begin{equation}
\mathcal{L}_2[0,1] =\left\lbrace f:\int_0^1f^2dt <\infty \right\rbrace
\end{equation}
with an inner product $\langle f,g\rangle=\int_0^1fgdt$. 

Consider a regression problem with observations modeled as $y_i = f(t_i)+\varepsilon_i$, $i=1,\ldots,n$, where $\varepsilon_i\sim N(0,\sigma^2)$ are \iid Gaussian noise and $f\in\mathcal{C}^{(m)}[0,1]=\{f:f^{(m)}\in \mathit{L}_2[0,1]\}$. The classic nonparametric or semi-parametric regression is a function that minimizes the following penalized sum of squared functional 
\begin{equation}\label{GaussianProcessGeneralObjective}
\frac{1}{n}\sum_{i=1}^{n}\left( y_i-f(t_i) \right)^2 + \lambda \int_{0}^{1} \left( f^{(m)}\right)^2dt, 
\end{equation}
where the first term is the lack of fit of $f$ to the data. The parameter $\lambda$ in the second term is a fixed smoothing parameter controlling the trade-off between over-fitting and bias \citep{esl2009}. The minimizer $f_\lambda$ of the above equation resides in an $n$-dimensional space and the computation in multivariate settings is generally of the order $\BigO{n^3}$ \citep{kim2004smoothing}. \cite{schoenberg1964spline} shows that a piecewise polynomial smoothing spline of degree $2m-1$ provides an aesthetically satisfying method for estimating $f$ if $\mathbf{y}=\left\lbrace y_1,\ldots,y_n\right\rbrace$ cannot be interpolated exactly by some polynomial of degree less than $m$. For instance, when $m=2$, a piecewise cubic smoothing spline provides a powerful tool to estimate the above nonparametric function, in which the penalty term is $\int f''^2dt$ \citep{hastie1990generalized}.

Further, \cite{kimeldorf1971some, kimeldorf1970correspondence} explore the corresponding between smoothing spline and Bayesian estimation. Actually, \cite{wahba1978improper} shows that a Bayesian version of this problem is to take a Gaussian process prior $f(t_i) = a_0+a_1t_i+\cdots + a_{m-1}t_i^{m-1} + x_i$ on $f$ with $x_i=X(t_i)$ being a zero-mean Gaussian process whose $m$th derivative is scaled white noise, $i=1,\ldots,n$ \citep{speckman2003fully}. The extended Bayes estimates $f_\lambda$ with a ``partially diffuse'' prior is exactly the same as the spline solution. \cite{heckman1991minimax} show that if prior distribution of the vector  $\mathbf{f}=(f(t_1),\ldots,f(t_n))^\top$ is unknown but lies in a known class $\Omega$,
the estimator $\hat{f}$ is found by minimizing the $\max\E\lbrack \hat{f}-f\rbrack^2$. \cite{branson2017nonparametric} propose a Gaussian process regression method that acts as a Bayesian analog to local linear regression for sharp regression discontinuity designs. It is no doubt that one of the attractive features of the Bayesian approach is that, in principle, one can solve virtually any statistical decision or inference problem. Particularly, one can provide an accuracy assessment for $\hat{f}=\E (f\mid \mathbf{y})$ using posterior probability regions \citep{cox1993analysis}.

The problem of choosing the smoothing parameter is ubiquitous in curve estimation, and there are two different philosophical approaches to this question. The first one is to regard the free choice of smoothing parameter as an advantageous feature of the procedure. The other one is to find the parameter automatically by the data \citep{green1993nonparametric}. We prefer the latter one, to use data to train our model and to find the best parameters. The most well-known method is cross-validation.

Assuming that mean of the random errors is zero, the true regression curve $f(t)$ has the property that, if an observation $y$ is taken away at a point $t$, the value $f(t)$ is the best predictor of $y$ in terms of returning a least value of $\left(y-f(t)\right)^2$. 

Now, focus on an observation $y_i$ at point $t_i$ as being a new observation by omitting it from the set of data, which are used to estimate $\hat{f}$. Denote by $\hat{f}^{(-i)}(t,\lambda)$ the estimated function from the remaining data, where $\lambda$ is the smoothing parameter. Then $\hat{f}^{(-i)}\left(t,\lambda\right)$ is the minimizer of  
\begin{equation}\label{originalobjective}
\frac{1}{n}\sum_{j \neq i}\left(y_j-f(t_j) \right)^2+\lambda\int (f'')^2dt,
\end{equation}
and can be quantified by the cross-validation score function
\begin{equation}
\mbox{CV}(\lambda)=\frac{1}{n}\sum_{i=1}^{n}\left(  y_i-\hat{f}^{(-i)}(t_i,\lambda)\right) ^2.
\end{equation}
The basis idea of the cross-validation is to choose the value of $\lambda$ that minimizes $\mbox{CV}(\lambda)$ \citep{green1993nonparametric}. 

An efficient way to calculate the cross-validation score is introduced by \cite{green1993nonparametric}. Because of the value of the smoothing spline $\hat{f}$ depending linearly on the data $y_i$, define the matrix $A(\lambda)$, which is a map vector of observed values $y_i$ to predicted value $\hat{f}(t_i)$. Then we have
\begin{equation}\label{crossvalidationmatrixA}
\hat{\mathbf{f}}=A(\lambda)\mathbf{y}.
\end{equation}
Based on the correspondence between nonparametric regression and Bayesian estimation, \cite{craven1978smoothing} propose a generalized cross-validation estimate for the minimizer $f_\lambda$. The estimated $\hat{\lambda}$ is the minimizer of the function where the trace of matrix $A(\lambda)$ in \eqref{crossvalidationmatrixA} is incorporated. It is also possible to establish an optimal convergence property for the estimator when the number of observations in a fixed interval tends to infinity \citep{wecker1983signal}. A highly efficient algorithm to optimize generalized cross-validation and generalized maximum likelihood scores with multiple smoothing parameters via the Newton method was proposed by \cite{gu1991minimizing}. This algorithm can also be applied to maximum likelihood estimation and restricted maximum likelihood estimation. The behavior of the optimal regularization parameter in different regularization methods was investigated by \cite{wahba1990optimal}.

In this paper, we prove that the V-spline, which incorporates both $f$ and $f'$ but penalizes excessive $f''$ in the penalty term, can be estimated by a Bayesian approach in a certain reproducing kernel Hilbert space. An extended GCV is used to find the optimal parameters for the V-spline.

\section{Polynomial Smoothing Splines on $[0, 1]$ as Bayes Estimates}

A polynomial smoothing spline of degree $2m-1$ is a piecewise polynomial of the same degree on each interval $[t_i,t_{i+1})$, $i=1, \ldots, n-1$, and the first $2m-2$ derivatives are continuous at the joined points. For instance, when $m=2$,  a piecewise cubic smoothing spline is a special case of the polynomial smoothing spline providing a powerful tool to estimate the objective function \eqref{GaussianProcessGeneralObjective} in the space  $\mathcal{C}^{(2)}[0,1]$, where the penalty term is $\int f''^2dt$ \citep{hastie1990generalized, wang1998smoothing}. If a general space $\mathcal{C}^{(m)}[0,1]$ is equipped with an appropriate inner product, it can be made as a reproducing kernel Hilbert space. 

\subsection{Polynomial Smoothing Spline}

A spline is a numeric function that is piecewise-defined by polynomial functions, which possesses a high degree of smoothness at the places where the polynomial pieces connect (known as \textit{knots}) \citep{judd1998numerical, chen2009feedback}. Suppose we are given a set of paired data $(t_1,y_1),\ldots, (t_n,y_n)$ on the interval $[0,1]$, satisfying $0< t_1 < \cdots <t_n < 1$. A piecewise polynomial function $f(t)$ can be obtained by dividing the interval into contiguous intervals $(t_1,t_2),\ldots,(t_{n-1},t_n)$ and represented by a separate polynomial on each interval. For any continuous $f\in \mathcal{C}^{(m)}[0,1]$, it can be represented in a linear combination of basis functions $h_m(t)$ as $f(t) =\sum_{m=1}^{M}\beta_mh_m(t)$, where $\beta_m$ are coefficients \citep{ellis2009}. It is just like every vector in a vector space can be represented as a linear combination of basis vectors.

A smoothing polynomial spline is uniquely the smoothest function that achieves a given degree of fidelity to a particular data set \citep{whittaker1922new}. In deed, the minimizer of function \eqref{GaussianProcessGeneralObjective} is the curve estimate $\hat{f}(t)$ over all spline functions $f(t)$ with $m-1$ continuous derivatives fitting observed data in the space $\mathcal{C}^{(m)}[0,1]$. In fact, the representer theorem \citep{kimeldorf1971some} tells us that the function $f_{\min} =\arg \min \frac{1}{n}\sum_i\left(y_i-f(t_i)\right)^2+\lambda \Vert f\Vert_C^2$ can be represented in the form $f_{\min} =\sum_i\alpha_iK(\cdot,t_i)$, a linear combination of a  positive-definite real-valued kernel $K(\cdot,\cdot)$ at each point. This is true for any arbitrary loss function \citep{rudin2005stability}.

Further, \cite{wahba1978improper} proves that if $f(t)$ has the prior distribution which is the same as the distribution of the stochastic process $\chi(t)$ on $[0,1]$, 
\begin{equation}\label{gausspriorequation}
\chi(t)=\sum_{\nu=1}^m d_\nu \phi_\nu(t)+b^{\frac{1}{2}}Z(t),
\end{equation}
where $d_i\sim N(0,\xi I)$, $Z(t)=\int_0^t\frac{(t-u)^{m-1}}{(m-1)!}dW(u)$ is the integrated Wiener process, then the polynomial spline $f_\lambda$ is the minimizer of the objective function \eqref{GaussianProcessGeneralObjective} having the property that 
\begin{equation}
f_\lambda(t)=\lim\limits_{\xi\rightarrow\infty} \E_\xi \left\lbrace f(t)\mid \mathbf{f}=\mathbf{y}  \right\rbrace,
\end{equation}
with $\lambda=\sigma^2/nb$, where $\E_\xi$ is expectation over the posterior distribution of $f(t)$ with the prior \eqref{gausspriorequation}. $\xi=\infty$ corresponds to the ``diffuse'' prior on $d$.

\subsection{Reproducing Kernel Hilbert Space in $C^{(m)}[0,1]$}

Any $f\in \mathcal{C}^{(m)}[0,1]$ has a standard Taylor expansion, which is 
\begin{equation}
f(t) = \sum_{\nu=0}^{m-1}\frac{t^\nu}{\nu!}f^{(\nu)}(0) + \int_{0}^{1}\frac{(t-u)_+^{m-1}}{(m-1)!}f^{(m)}(u)du,
\end{equation}
where $(\cdot)_+ =\max\lbrace0, \cdot\rbrace$. With an inner product 
\begin{equation}
\langle f,g \rangle = \sum_{\nu=0}^{m-1}f^{(\nu)}(0)g^{(\nu)}(0) +  \int_{0}^{1}f^{(m)}(t) g^{(m)}(t)dt,
\end{equation}
the representer is 
\begin{equation}\label{GaussianProcessKernelR}
\begin{split}
R_s(t) &=\sum_{\nu=0}^{m-1} \frac{s^{\nu}}{\nu!} \frac{t^{\nu}}{\nu!} +\int_0^1\frac{ (s-u)_+^{m-1}}{(m-1)!} \frac{ (t-u)_+^{m-1}}{(m-1)!} du \\ 
&= R_0(s,t)+R_1(s,t). 
\end{split}
\end{equation}
It is easy to prove that $R(s,t)$ is a non-negative reproducing kernel, by which $\langle R(s,t),f(t) \rangle = \langle R_s(t),f(t) \rangle=f(s)$. Additionally, $R_s^{(\nu)}(0) = s^\nu/\nu!$ for $\nu = 0,\ldots, m-1$.

Before moving on to further steps, we are now introducing the following two theorems. 
\begin{theorem}\citep{aronszajn1950theory}\label{theoremRKHS}
Suppose $R$ is a symmetric, positive definite kernel on a set $X$. Then, there is a unique Hilbert space of functions on $X$ for which $R$ is a reproducing kernel. 
\end{theorem}
\begin{theorem}\citep{gu2013smoothing}\label{theoremKernel}
If the reproducing kernel $R$ of a space $\mathcal{H}$ on domain $X$ can
be decomposed into $R = R_0 + R_1$, where $R_0$ and $R_1$ are both non-negative definite, $R_0(x, \cdot),R1(x,\cdot) \\ \in \mathcal{H}$, for $ \forall x \in X$, and $\langle R_0(x, \cdot),R_1(y, \cdot) \rangle= 0$, for $\forall x, y \in X$, then the spaces $\mathcal{H}_0$ and $\mathcal{H}_1$ corresponding respectively to $R_0$ and $R_1$ form a tensor sum decomposition of $\mathcal{H}$. Conversely, if $R_0$ and $R_1$ are both  nonnegative definite and $\mathcal{H}_0 \cap \mathcal{H}_1 =\left\lbrace 0\right\rbrace$, then $\mathcal{H} =\mathcal{H}_0 \oplus \mathcal{H}_1$ has a reproducing kernel $R = R_0 + R_1$.
\end{theorem}

According to Theorem \ref{theoremRKHS}, the Hilbert space associated with $R(\cdot)$ can be constructed as containing all finite linear combinations of the form $\sum a_iR(t_i,\cdot)$, and their limits under the norm induced by the inner product $\langle R(s,\cdot),R(t,\cdot) \rangle = R(s,t)$. As for Theorem \ref{theoremKernel}, it is easy to verify that $R_0$ corresponds to the space of polynomials $\mathcal{H}_0 =\left\lbrace f:f^{(m)}=0\right\rbrace$ with an inner product $\langle f, g \rangle_0 = \sum_{\nu=0}^{m-1} f^{(\nu)}(0)g^{(\nu)}(0)$ and $R_1$ corresponds to the orthogonal complement of $\mathcal{H}_0$, that is $\mathcal{H}_1 =\left\lbrace   f:f^{(\nu)}(0)=0,\nu = 0, \ldots,m-1, \int_{0}^{1}(f^{(m)})^2dt <\infty  \right\rbrace$ with an inner product $\langle f, g\rangle_1 = \int_{0}^{1}f^{(m)}g^{(m)}dt $. 

Given a set of sampling points, any $f\in C^{(m)}[0,1]$ has the following form 
\begin{equation}
f(t)=\sum_{\nu=1}^m d_\nu \phi_\nu(t)+\sum_{i=1}^n c_iR_1(t,t_i).
\end{equation}
where $\left\lbrace \phi_\nu(t)\right\rbrace$ is a set of basis functions of space $\mathcal{H}_0$ and $R(\cdot,t)$ is the representer in $\mathcal{H}_1$ \citep{wang2011smoothing}. 

Additionally, the coefficients $c_i$ and $d_\nu$ might be changed when different $\phi_\nu$ and $R_1$ are used, but the function estimate remains the same regardless of the choices of $\phi_\nu$ and $R_1$ \citep{gu2013smoothing}.

\subsection{Polynomial Smoothing Splines as Bayes Estimates}

Because it is possible to interpret the smoothing spline regression estimator as a Bayes estimate when the mean function $r(\cdot)$ is given an improper prior distribution \citep{wahba1990spline, berlinet2011reproducing}. Therefore, one can find that the posterior mean of $f$ on $[0,1]$ with a vague improper prior is the polynomial smoothing spline of the objective function  \eqref{GaussianProcessGeneralObjective}. 

Consider $f=f_0+f_1$ on $[0,1]$, with $f_0$ and $f_1$ having independent Gaussian priors with zero means and covariances satisfying  
\begin{align}
\E \lbrack f_0(s)f_0(t)\rbrack  &= \tau^2 R_0(s,t)=\tau^2 \sum_{\nu=0}^{m-1}\frac{s^\nu}{\nu!}\frac{t^\nu}{\nu!},\label{Ef0f0} \\
\E \lbrack f_1(s)f_1(t) \rbrack &= bR_1(s,t) = b\int_{0}^{1} \frac{(s-u)_+^{m-1}}{(m-1)!} \frac{(t-u)_+^{m-1}}{(m-1)!},\label{Ef1f1}
\end{align}
where $R_0$ and $R_1$ are from \eqref{GaussianProcessKernelR}. Because of the observations are normally distributed as $y_i\sim N(f(t_i),\sigma^2)$, then the joint distribution for $\mathbf{y}$ and $f(t)$ is normal with mean of zero and covariance matrix of 
\begin{align*}\Cov (f,\mathbf{y}) = 
\begin{bmatrix}
bQ+\tau SS^\top+\sigma^2 I & b\xi +\tau^2 S\phi \\
b\xi^\top + \tau^2\phi^\top S^\top & bR_1(t,t) +\tau^2\phi^\top \phi
\end{bmatrix},
\end{align*}
where $\left\lbrace Q_{i,j}\right\rbrace_{n\times n}=R_1(t_i,t_j)$, $\left\lbrace S_{i,\nu}\right\rbrace_{n\times m}=t_i^{\nu-1}/(\nu-1)!$, $\left\lbrace \xi_{i,1}\right\rbrace_{n\times 1}=R_1(t_i,t)$ and $\left\lbrace \phi_{\nu,1}\right\rbrace_{m\times 1}=t^{\nu-1}/(\nu-1)!$. 
Consequently, the posterior is 
\begin{equation}
\begin{split}
\E \lbrack f(t)\mid\mathbf{y}\rbrack &= \left(b\xi^\top +\tau \phi^\top s^\top\right)\left(bQ+\tau^2 SS^\top+\sigma^2I\right)^{-1}\mathbf{y} \\
&= \xi^\top\left(Q+\rho SS^\top+n\lambda I\right)^{-1}\mathbf{y}+ \phi^\top\rho S^\top\left(Q+\rho SS^\top +n\lambda I\right)^{-1}\mathbf{y},
\end{split}
\end{equation}
where $\rho = \tau^2/b$ and $n\lambda=\sigma^2/b$. Furthermore, by denoting $M=Q+n\lambda I$, \cite{gu2013smoothing} gives that, when $\rho\rightarrow \infty$, the posterior mean is in the form $\E\lbrack f(t)\mid y_{1:n}\rbrack = \xi^\top\mathbf{c}+\phi^\top\mathbf{d}$ with coefficient vectors 
\begin{align}
\mathbf{c}&=\left(M^{-1}-M^{-1}S\left(S^\top M^{-1}S\right)^{-1}S^\top M^{-1}\right)\mathbf{y},\\
\mathbf{d}&=\left(S^\top M^{-1}S\right)^{-1}S^\top M^{-1}\mathbf{y}.
\end{align}

\begin{theorem}\citep{gu2013smoothing}
The polynomial smoothing spline of \eqref{GaussianProcessGeneralObjective} is the posterior mean of $f = f_0 +f_1$, where $f_0$ diffuses in span $\left\lbrace t^{\nu-1}, \nu= 1, \ldots , m\right\rbrace$ and $f_1$ has a Gaussian process prior with mean zero and a covariance function
\begin{equation}
bR_1(s,t) = b\int_{0}^{1} \frac{\left(s-u\right)_+^{m-1}}{(m-1)!} \frac{\left(t-u\right)_+^{m-1}}{(m-1)!},
\end{equation}
for $b=\sigma^2/n\lambda$. 
\end{theorem}

\textit{Remark}: Equation \eqref{Ef0f0} can be obtained from equation \eqref{gausspriorequation} if we assume $d_\nu \sim N\left(0,\tau^2I_{m\times m}\right)$. Therefore the limit of $\rho=\tau^2/b\to\infty$ indicates a diffuse prior for the coefficients $\mathbf{d}$.

\subsection{Gaussian Process Regression}

Gaussian processes are the extension of multivariate Gaussian to infinite-sized collections of real value variables, any finite number of which have a joint Gaussian distribution \citep{rasmussen2006gaussian}. Gaussian process regression is a probability distribution over functions. It is fully defined by its mean $m(t)$ and covariance $K(s,t)$ function as 
\begin{align}
m(t)&=\E \lbrack f(t)\rbrack \\
K(s,t)&=\E \lbrack \left(f(s)-m(s)\right) \left(f(t)-m(t)\right)\rbrack,
\end{align}
where $s$ and $t$ are two variables. A function $f$ distributed as such is denoted in form of 
\begin{equation}
f \sim GP\left(m(t),K(s,t) \right).
\end{equation}
Usually the mean function is assumed to be zero everywhere. 

Given a set of input variables $\mathbf{t} = \left\lbrace t_1,\ldots,t_n\right\rbrace$ for function $f(t)$ and the output $\mathbf{y}=f(\mathbf{t})+\varepsilon$ with \iid  Gaussian noise $\varepsilon$ of variance $\sigma_n^2$,  we can use the above definition to predict the value of the function $f_*=f(t_*)$ at a particular input $t_*$. As the noisy observations becoming
\begin{equation} 
\Cov (y_p,y_q) = K(t_p,t_q)+\sigma_n^2 \delta_{pq}
\end{equation}
where $\delta_{pq}$ is a Kronecker delta which is one if and only if $p=q$ and zero otherwise, the joint distribution of the observed outputs $\mathbf{y}$ and the estimated output $f_*$ according to prior is
\begin{equation}
\begin{bmatrix}
\mathbf{y}\\
f_*
\end{bmatrix} \sim N \left(  
0,  \begin{bmatrix}
K(\mathbf{t},\mathbf{t}) +\sigma_n^2I& K(\mathbf{t},t_*) \\
K(t_*,\mathbf{t}) & K(t_*,t_*)
\end{bmatrix} 
\right).
\end{equation}
The posterior distribution over the predicted value is obtained by conditioning on the observed data
\begin{equation}
f_* \mid  \mathbf{y},\mathbf{t},t_* \sim N\left(\bar{f_*},\Cov (f_*)\right)
\end{equation}
where 
\begin{align}
\bar{f_*}&=\E \left( f_* \mid  \mathbf{y},\mathbf{t},t_*\right) = K(t_*,\mathbf{t})\left( K(\mathbf{t},\mathbf{t})+\sigma_n^2\right) ^{-1}\mathbf{y},\\
\Cov(f_*)&=K(t_*,t_*)-K(t_*,\mathbf{t})\left( K(\mathbf{t},\mathbf{t})+\sigma_n^2I\right) ^{-1}K(\mathbf{t},t_*).
\end{align}
Therefore it can seen that the Bayesian estimation of a smoothing spline is a special format of Gaussian process regression with diffuse prior and the covariance matrix $R(s,t)$.

\section{V-Splines and Bayes Estimate}

\subsection{V-Splines}

In a nonparametric regression, consider $n$ paired time series points $\left\lbrace t_1,y_1,v_1\right\rbrace$, $\ldots$, $\left\lbrace t_n,y_n,v_n\right\rbrace$, such that $0 < t_1< \cdots < t_n < 1$, $y$ is the position information and $v$ indicates its velocity. As in \citep{silverman1985some} and \citep{donoho1995wavelet}, we use a positive penalty function $\lambda(t)$ in the following objective function rather than a constant $\lambda$ in \eqref{originalobjective}. 

Given function $f:[0,1]\mapsto \mathbb{R}$ and $\gamma>0$, define the objective function 
\begin{equation}\label{tractorsplineObjective}
J[f]= \frac{1}{n} \sum_{i=1}^{n} \left( f(t_i)-y_i \right)^2 + \frac{\gamma}{n} \sum_{i=1}^{n} \left( f'(t_i)-v_i \right)^2 +\sum_{i=1}^{n-1} \int_{t_i}^{t_{i+1}}\lambda(t)  f''^2(t)dt,
\end{equation}
where $\gamma$ is the parameter that weights the residuals between $\mathbf{f}'=\{f'(t_1),\ldots,f'(t_n)\}$ and $\mathbf{v}=\{v_1,\ldots,v_n\}$. We make a  simple assumption that $\lambda(t)$ is a piecewise constant and adopts a constant value $\lambda_i$ on interval $(t_i,t_{i+1})$ for $i=1,\ldots, n-1$. 
\begin{theorem}\label{TractorSplineTheorem}
For $n\geq2$, the objective function $J[f]$ is minimized by a cubic spline that is unique and linear outside the knots.
\end{theorem}
A further minimizer of \eqref{tractorsplineObjective} is named \textit{V-spline}, coming from the incorporation with velocity information and applications on vehicle and vessel tracking. It is the solution to the objective function \eqref{tractorsplineObjective}, where an extra term for $f'(t)-v$ and an extra parameter $\gamma$ are incorporated. The penalty parameter $\lambda(t)$ is a function varying on different domains. If $\lambda(t)$ is constant and $\gamma=0$, the V-spline degenerates to a conventional cubic smoothing spline consisting of a set of given basis functions.

However, the Bayes estimate for a polynomial smoothing spline requires a constant penalty parameter. For this constraint, it is assumed that $\lambda(t)$ stays the same on each subinterval in $[0,1]$ and named the solution ``trivial V-spline''.  In this section, we still use ``V-spline'' for sake of simplicity.

\subsection{Reproducing Kernel Hilbert Space $\mathcal{C}_{\mbox{\scriptsize p.w.}}^{(2)}[0,1]$}

The space $\mathcal{C}^{(m)}[0,1]=\left\lbrace  f:f^{(m)}\in \mathit{L}_2[0,1] \right\rbrace$ is a set of functions $f$ whose $m$th derivatives are square integrable on the domain $[0,1]$. For a V-spline, it only requires $m=2$. In fact, its second derivative is piecewise linear but is not necessarily continuous at the knots. Besides, if and only if $\lambda(t)$ is constant and $\gamma=0$, the second derivative is piecewise linear and continuous at the knots. Here we are introducing the space 
\begin{equation*}
\mathcal{C}_{\mbox{\scriptsize p.w.}}^{(2)}[0,1]=\left\lbrace f: f''\in \mathit{L}_2[0,1], f,f' \mbox{ are continuous and } f'' \mbox{ is piecewise linear}\right\rbrace,
\end{equation*}
in which the second derivative of any function $f$ is not necessarily continuous.

Given a sequence of paired data $\left\lbrace (t_1,y_1,v_1),\ldots, (t_n,y_n,v_n) \right\rbrace$, the the minimizer of 
\begin{equation}\label{maineq}
J[f]=\frac{1}{n}\sum_{i=1}^{n}(y_i-f(t_i))^2+\frac{\gamma}{n}\sum_{i=1}^{n}(v_i-f'(t_i))^2+\lambda \int_{0}^{1}f''^2dt
\end{equation}
in the space $\mathcal{C}_{\mbox{\scriptsize p.w.}}^{(2)}[0,1]$ is a V-spline. Equipped with an appropriate inner product
\begin{equation}\label{TractorSplineInnerProduct}
\langle f,g \rangle=f(0) g(0)+f'(0) g'(0)+\int_{0}^{1}f''(t)g''(t)dt,
\end{equation}
the space $\mathcal{C}_{\mbox{\scriptsize p.w.}}^{(2)}[0,1]$ is made a reproducing kernel Hilbert space. In fact, the representer $R_s(\cdot)$ is 
\begin{equation}\label{kerneleq}
R_s(t)=1+st+\int_{0}^{1} (s-u)_+(t-u)_+du.
\end{equation}
It can be seen that $R_s(0)=1, R'_s(0)=s$, and $R''_s(t)=(s-t)_+$. The two terms of the reproducing kernel $R(s,t)=R_s(t)\triangleq R_0(s,t)+R_1(s,t)$, where
\begin{align} \label{TractorSplineKernelR0}
R_0(s,t)&=1+st \\ \label{TractorSplineKernelR1}
R_1(s,t)&=\int_{0}^{1} (s-u)_+(t-u)_+du
\end{align}
are both non-negative definite themselves.

According to Theorem \ref{theoremKernel}, $R_0$ can correspond the space of polynomials $\mathcal{H}_0=\left\lbrace f:f''=0\right\rbrace$ with an inner product $\langle f,g \rangle_0= f(0)g(0)+f'(0)g'(0)$, and $R_1$ corresponds the orthogonal complement of $\mathcal{H}_0$
\begin{equation}
\mathcal{H}_1=\left\lbrace f:f(0)=0, f'(0)=0, \int_{0}^{1}f''(t)^2dt<\infty\right\rbrace
\end{equation}
with inner product $\langle f,g \rangle_1=\int_{0}^{1}f''g''dt$. Thus, $\mathcal{H}_0$ and $\mathcal{H}_1$ are two subspaces of the $\mathcal{C}_{\mbox{\scriptsize p.w.}}^{(2)}[0,1]$, and the reproducing kernel is $R_s(\cdot) = R_0(s,\cdot)+R_1(s,\cdot)$.

Define a new notation $\dot{R}(s,t)=\frac{\partial R}{\partial s}(s,t)=\frac{\partial R_0}{\partial s}(s,t)+\frac{\partial R_1}{\partial s}(s,t)=t+\int_0^s(t-u)_+du$. Obviously $\dot{R}_s(t) \in \mathcal{C}_{\mbox{\scriptsize p.w.}}^{(2)}[0,1]$. Additionally, we have $\dot{R}_s(0)=0, \dot{R}'_s(0)=\frac{\partial \dot{R}_s}{\partial t}(0)=1$, and $ \dot{R}''_s(t)=\begin{cases}
0 & s\leq t \\ 1 & s>t \end{cases}$. Then, for any $f\in \mathcal{C}_{\mbox{\scriptsize p.w.}}^{(2)}[0,1]$, we have 
\begin{equation}
\langle \dot{R}_s,f\rangle =\dot{R}_s(0)f(0)+\dot{R}'_s(0)f'(0)+\int_0^1\dot{R}''_s f''	(u) du=f'(0)+\int_0^t f''(u)du=f'(t).
\end{equation}
%
It can be seen that the first term $\dot{R}_0=t\in \mathcal{H}_0$, and the space spanned by the second term $\dot{R}_1=\int_0^s(t-u)_+du$, denoted as $\mathcal{\dot{H}}$, is a subspace of $\mathcal{H}_1$, and $\mathcal{\dot{H}} \ominus \mathcal{H}_1\neq \emptyset$. Given the sample points $t_j, j=1, \ldots, n$, in equation \eqref{maineq} and noting that the space
\begin{equation}
\mathcal{A}=\left\lbrace f: f=\sum_{j=1}^{n}\alpha_jR_1(t_j,\cdot)+\sum_{j=1}^{n}\beta_j\dot{R}_1(t_j,\cdot)\right\rbrace 
\end{equation}
is a closed linear subspace of $\mathcal{H}_1$. Then, we have a new space $\mathcal{H}_*=\mathcal{\dot{H}} \cup \mathcal{A}$. Thus, the two new sub spaces in $\mathcal{C}_{\mbox{\scriptsize p.w.}}^{(2)}[0,1]$ are $\mathcal{H}_0$ and $\mathcal{H}_*$.

For any $f\in\mathcal{C}_{\mbox{\scriptsize p.w.}}^{(2)}[0,1]$, it can be written as 
\begin{equation}\label{GaussianProcessFunctionF}
f(t)=d_1+d_2t+\sum_{j=1}^{n}c_jR_1(t_j,t)+\sum_{j=1}^{n}b_j\dot{R}_1(t_j,\cdot) +\rho(t)
\end{equation}
where $\mathbf{d}=\lbrace d_1,d_2\rbrace,\mathbf{c}=\lbrace c_j\rbrace$ and $\mathbf{b}=\lbrace b_j\rbrace$, $j=1,\ldots,n$, are coefficients, and $\rho(t) \in \mathcal{H}_1 \ominus \mathcal{H}_*$. Thus, by substituting to the equation \eqref{maineq}, it can be written as 
\begin{equation}\label{GassianProcessRawequation}
\begin{split}
nJ[f]=&\sum_{i=1}^n \left( y_i - d_1-d_2t-\sum_{j=1}^{n}c_jR_1(t_j,t_i)-\sum_{j=1}^{n}b_j\dot{R}_1(t_j,t_i)-\rho(t_i) \right) ^2\\
+&\gamma\sum_{i=1}^n \left( v_i - d_2-\sum_{j=1}^{n}c_jR'_1(t_j,t_i)-\sum_{j=1}^{n}b_j\dot{R}'_1(t_j,t_i)-\rho'(t_i) \right) ^2\\
+&n\lambda \int_0^1 \left( \sum_{j=1}^{n}c_jR''_1(t_j,t)+\sum_{j=1}^{n}b_j\dot{R}''_1(t_j,t)+\rho''(t)\right)^2dt
\end{split}
\end{equation}
Because of orthogonality, $\rho(t_i) = \langle R_1(t_i,\cdot),\rho\rangle=0$, $\rho'(t_i) = \langle\dot{R}_1(t_i,\cdot),\rho'\rangle=0$, $i=1,\ldots,n$. By denoting that 
\begin{align*}
S&=\left\lbrace S_{ij} \right\rbrace_{n\times 2}=\begin{bmatrix}1 & t_i \end{bmatrix} ,& Q&=\left\lbrace Q_{ij} \right\rbrace_{n\times n}= R_1(t_j,t_i), & P&=\left\lbrace P_{ij} \right\rbrace_{n\times n}= \dot{R}_1(t_j,t_i), \\
S'&=\left\lbrace S'_{ij} \right\rbrace_{n\times 2}=\begin{bmatrix} 0 & 1 \end{bmatrix} ,& Q'&=\left\lbrace Q'_{ij} \right\rbrace_{n\times n}= R_1'(t_j,t_i), & P'&=\left\lbrace P'_{ij} \right\rbrace_{n\times n}= \dot{R}_1'(t_j,t_i). 
\end{align*}
and noting that $\int_0^1R''_1(t_i,t)R''_1(t_j,t)dt=R_1(t_i,t_j)$, $\int_0^1R''_1(t_i,t)\dot{R}''_1(t_j,t)dt=\int_0^{v}(t_i-t)dt=\dot{R}_1(t_j,t_i)$, and $\int_0^1\dot{R}''_1(t_i,t)\dot{R}''_1(t_j,t)dt=\int_0^{v}1dt=\dot{R}'_1(t_i,t_j)$, where $v=\mbox{min}\lbrace t_i,t_j\rbrace$, the above equation \eqref{GassianProcessRawequation} can be written as 
\begin{equation}\label{matriteq}
\begin{split}
nJ[f]=&\left(\mathbf{y}-S\mathbf{d}-Q\mathbf{c}-P\mathbf{b}\right)^\top \left(\mathbf{y}-S\mathbf{d}-Q\mathbf{c}-P\mathbf{b}\right)\\
+&\gamma\left(\mathbf{v}-S'\mathbf{d}-Q'\mathbf{c}-P'\mathbf{b}\right)^\top \left(\mathbf{v}-S'\mathbf{d}-Q'\mathbf{c}-P'\mathbf{b}\right)\\
+&n\lambda \left(\mathbf{c}^\top Q\mathbf{c} + 2\mathbf{c}^\top P\mathbf{b}+ \mathbf{b}^\top P'\mathbf{b}\right)+n\lambda\left(\rho,\rho\right).
\end{split}
\end{equation}
Note that $\rho$ only appears in the third term and is minimized at $\rho=0$. Hence, a V-spline resides in the space $\mathcal{H}_0\oplus \mathcal{H}_*$ of finite dimension. Thus, the solution to \eqref{maineq} is computed via the minimization of the first three terms in \eqref{matriteq} with respect to $\mathbf{d}$, $\mathbf{c}$ and $\mathbf{b}$.

\subsection{Posterior of Bayes Estimates}\label{sectionBayesEstimate}


In a general process, we know that $p(\mathbf{y},\mathbf{v}\mid f) = N(f,\Gamma)$, where $\Gamma$ is a covariance matrix. However, we are more interested in $f$ given measurements, which is  
\begin{equation}
p(f\mid \mathbf{y},\mathbf{v}) \propto p(\mathbf{y},\mathbf{v}\mid f)p(f),
\end{equation}
where $f\sim GP(0,\Sigma)$ is a Gaussian process prior. In fact, the covariance matrix $\Sigma$ is associated to the inner product $R(s,t)$. 

Observing $y_i\sim N\left(f (t_i),\sigma^2\right)$ and $v_i\sim N\left(f (t_i),\frac{\sigma^2}{\gamma}\right)$, $i=1,\ldots,n$, the joint distribution of $\mathbf{y},\mathbf{v}$ and $f(t)$ is normal with mean zero and a covariance matrix can be found by the following 
\small
\begin{equation}
\begin{aligned}
\E\lbrack f(s)f (t)\rbrack &=\tau^2R_0(s,t)+\beta R_1(s,t) & \E \lbrack f(s)f '(t)\rbrack&=\tau^2R_0'(s,t)+\beta R_1'(s,t) \\
\E\lbrack f '(s)f (t)\rbrack&=\tau^2\dot{R}_0(s,t)+\beta\dot{R}_1(s,t) & \E\lbrack f '(s)f '(t)\rbrack&=\tau^2\dot{R}'_0(s,t)+\beta\dot{R}'_1(s,t) \\
\footnotesize \E\lbrack y_i,y_j\rbrack=\tau^2 & R_0(s_i,s_j)+\beta R_1(s_i,s_j)+\sigma^2\delta_{ij}   & 
\E \lbrack v_i,v_j\rbrack=\tau^2&\dot{R}'_0(s_i,s_j)+\beta\dot{R}'_1(s_i,s_j) +\frac{\sigma^2}{\gamma}\delta_{ij} \\ 
\normalsize
\E\lbrack v_i,y_j\rbrack&=\tau^2\dot{R}_0(s_i,s_j)+\beta \dot{R}_1(s_i,s_j) &
\E\lbrack y_i,v_j\rbrack&=\tau^2R_0'(s_i,s_j)+\beta R_1'(s_i,s_j)\\
\E\lbrack y_i,f(s)\rbrack&=\tau^2 R_0(s_i,s)+\beta R_1(s_i,s)  & \E\lbrack y_i,f '(s)\rbrack&=\tau^2 R'_0(s_i,s)+\beta R'_1(s_i,s)  \\
\E\lbrack v_i,f(s)\rbrack&=\tau^2 \dot{R}_0(s_i,s)+\beta \dot{R}_1(s_i,s) & \E\lbrack v_i,f '(s)\rbrack&=\tau^2\dot{R}'_0(s_i,s)+\beta \dot{R}'_1(s_i,s)
\end{aligned}
\end{equation}
\normalsize where $R_0(s,t)$ and $R_1(s,t)$ are taken from \eqref{TractorSplineKernelR0} and \eqref{TractorSplineKernelR1}. 

Therefore, by using a standard result on multivariate normal distribution (such as Result 4.6 in \citep{johnson1992applied}), the posterior mean of $f(t)$ is seen to be 
\begin{align}\label{GPrhoeq}\footnotesize
\begin{split}
\E \lbrack f \mid  \mathbf{\mathbf{y}},\mathbf{v}\rbrack & =
\begin{bmatrix}
\Cov (\mathbf{y},f) & \Cov (f,\mathbf{\mathbf{v}})
\end{bmatrix}\begin{bmatrix}
\Var(\mathbf{y}) & \Cov(\mathbf{y},\mathbf{v})\\
\Cov(\mathbf{v},\mathbf{y}) & \Var(\mathbf{v})
\end{bmatrix}^{-1}\begin{bmatrix}
\mathbf{y}\\\mathbf{v}
\end{bmatrix}\\
&=
\begin{bmatrix}
\tau^2 \phi^\top S^\top+\beta \xi^\top & \tau^2  \phi^\top S'^\top+\beta \psi^\top 
\end{bmatrix}\begin{bmatrix}
\tau^2 SS^\top+\beta Q+\sigma^2 I& \tau^2 SS'^\top+\beta P\\
\tau^2 S'S^\top+\beta Q'& \tau^2 S'S'^\top+\beta P'+\frac{\sigma^2}{\gamma}I
\end{bmatrix}^{-1}\begin{bmatrix}
\mathbf{y}\\\mathbf{v}
\end{bmatrix}\\
&=
\begin{bmatrix}
\rho\phi^\top S^\top+ \xi^\top & \rho\phi^\top S'^\top+\psi^\top
\end{bmatrix}\begin{bmatrix}
\rho SS^\top+Q+n\lambda I& \rho SS'^\top+P\\
\rho S'S^\top+Q'& \rho S'S'^\top+P'+\frac{n\lambda}{\gamma}I
\end{bmatrix}^{-1}\begin{bmatrix}
\mathbf{y}\\\mathbf{v}
\end{bmatrix}\\
&=\left(\phi^\top \rho 
\begin{bmatrix} S\\ S' \end{bmatrix}^\top + \begin{bmatrix} \xi^\top & \psi^\top\end{bmatrix}\right)
\left(\rho\begin{bmatrix} S \\ S' \end{bmatrix}^\top \begin{bmatrix} S \\ S' \end{bmatrix}+
\begin{bmatrix} Q+n\lambda I& P\\
Q'& P'+\frac{n\lambda}{\gamma}I\end{bmatrix}\right) ^{-1}
\begin{bmatrix}\mathbf{y}\\ \mathbf{v} \end{bmatrix}\\
&\triangleq\phi^\top \rho T^\top \left(\rho T^\top T+M\right) ^{-1} \begin{bmatrix}\mathbf{y}\\ \mathbf{v} \end{bmatrix}
+ \begin{bmatrix} \xi^\top & \psi^\top\end{bmatrix}\left(\rho T^\top T+M\right) ^{-1} \begin{bmatrix}\mathbf{y}\\ \mathbf{v} \end{bmatrix}
\end{split}
\end{align}\normalsize
where $\phi$ is $2 \times 1$ matrix with entry $1$ and $t$, $\xi$ is $n\times 1$ matrix with $i$th entry $R(t_i,t)$, $T^\top =\begin{bmatrix} S^\top & S'^\top \end{bmatrix}$ and $\psi$ is $n\times 1$ matrix with $i$th entry  $\dot{R}(t_i,t)$, $\rho=\tau^2/\beta$ and $n\lambda =\sigma^2/\beta$. 

\begin{lemma}\label{GPLemma}
	Suppose $M$ is symmetric and nonsingular and $T$ is of full column rank. 
	\begin{align}
	&\lim\limits_{\rho \rightarrow \infty}\left(\rho TT^\top+M\right)^{-1}=M^{-1}-M^{-1}T\left(T^\top M^{-1}T\right)^{-1}T^\top M^{-1},\\
	&\lim\limits_{\rho \rightarrow \infty}\rho T^\top\left(\rho TT^\top+M\right)^{-1}=\left(T^\top M^{-1}T\right)^{-1}T^\top M^{-1}.
	\end{align}
\end{lemma}

Setting $\rho \rightarrow \infty$ in equation \eqref{GPrhoeq} and applying Lemma \ref{GPLemma}, the posterior mean $\E (f(t)\mid \mathbf{y},\mathbf{v})$ is $\hat{f}  = \mathbf{\phi}^\top \mathbf{d}+\mathbf{\xi}^\top \mathbf{c}+\mathbf{\psi}^\top \mathbf{b}$, with the coefficients given by
\begin{align} 
\mathbf{d}&=\left(T^\top M^{-1}T\right)^{-1}T^\top M^{-1}\begin{bmatrix}\mathbf{y} \\ \mathbf{v} \end{bmatrix},\\
\begin{bmatrix}\mathbf{c}\\ \mathbf{b}\end{bmatrix} &=
\left(M^{-1}-M^{-1}T\left(T^\top M^{-1} T\right)^{-1}T^\top M^{-1}\right)\begin{bmatrix}\mathbf{y}\\ \mathbf{v} \end{bmatrix},
\end{align} 
where $T=\begin{bmatrix} S\\S' \end{bmatrix}$ and $M=\begin{bmatrix} Q+n\lambda I& P\\
Q'& P'+\frac{n\lambda}{\gamma}I
\end{bmatrix}$.

It is easy to verify that $\mathbf{d},\mathbf{c},\mathbf{b}$ are the solutions to
\begin{equation}
\begin{cases}
S^\top \left(S\mathbf{d} +Q\mathbf{c}+P\mathbf{b}-\mathbf{y}\right) +\gamma S'^\top\left( S'\mathbf{d}+ P^\top \mathbf{c}+ P'\mathbf{b}-\mathbf{v}\right)=0, \\
Q\left(S\mathbf{d}+\left(Q+n\lambda I\right)\mathbf{c}+P\mathbf{b}-\mathbf{y}\right) + P \left( \gamma S' \mathbf{d} + \gamma P^\top \mathbf{c}+ \left(\gamma P'+n\lambda I\right) \mathbf{b}- \gamma \mathbf{v}\right)=0, \\
P^\top \left(S\mathbf{d}+\left(Q+n\lambda I\right) \mathbf{c} +P\mathbf{b}-\mathbf{y}\right)+P'\left(\gamma S'\mathbf{d}+P^\top \mathbf{c}+\left(\gamma P'+n\lambda I\right)\mathbf{b}- \gamma\mathbf{v}\right)=0.
\end{cases}
\end{equation}
Finally we obtain the following theorem: 
\begin{theorem}
The smoothing V-spline of \eqref{maineq} is the posterior
mean of $f=f_0+f_1 + \dot{f}_1$, where $f_0$ diffuses in span $\left\lbrace 1,t\right\rbrace$ and $f_1$, $\dot{f}_1$ have Gaussian process priors with mean zero and covariance functions
\begin{align}
\Cov \left(f_1,f_1\right)   &= \beta R_1\left(s,t\right)   =\beta \int_{0}^{1} \left(s-u\right)_+\left(t-u\right)_+du, \\
\Cov \left(\dot{f}_1,f_1\right)  &= \beta \dot{R}_1\left(s,t\right) =\beta \int_0^s\left(t-u\right)_+du,\\
\Cov \left(\dot{f}_1,\dot{f}_1\right)  &= \beta \dot{R}'_1\left(s,t\right) =\beta \min\lbrace s,t\rbrace,
\end{align}
for $\beta = \sigma^2/n\lambda$.
\end{theorem}

\section{Bayes Estimate for Non-trivial V-Spline}

For a sequence $0=t_0 < t_1<\cdots <t_n <t_{n+1}=1$ on the interval $[0,1]$ in the reproducing kernel Hilbert space  $\mathcal{C}_{\mbox{\scriptsize p.w.}}^{(2)}[0,1]$, define an inner product 
\begin{equation}\label{nontrivialInner}
\langle f,g\rangle = f(0)g(0)+f'(0)g'(0)+\sum_{i=0}^{n}w_i\int_{t_i}^{t_{i+1}}f''(t)g''(t)dt,
\end{equation}
where $w_i>0$, $i=0,\ldots,n$. The representer is 
\begin{align}
R_s(t) &= 1+st+\sum_{i=0}^{n}w_i^{-1}\int_{t_i}^{t_{i+1}}(s-u)_+(t-u)_+du, 
\end{align}
having the following properties 
\begin{align}
R'_s(t) &= s+\sum_{i=0}^{n}w_i^{-1}\int_{t_i}^{t_{i+1}}(s-u)_+\Theta(t-u) du,\\
\dot{R}_s(t) &= t+\sum_{i=0}^{n}w_i^{-1}\int_{t_i}^{t_{i+1}}\Theta(s-u)(t-u)_+ du,\\
R''_s(t) & = \sum_{i=0}^{n}w_i^{-1}\int_{t_i}^{t_{i+1}}(s-u)_+\delta(t-u)du, 
\end{align}
and $R_s(0)=1$, $R'_s(0)=s$. The function $\Theta(t-u)$ is the Heaviside function and $\delta(t-u)$ is the Dirac delta function. 

Further, $R(\cdot)$ and $\dot{R}(\cdot)$ on $\lbrack0,1\rbrack$ have the following properties 
\begin{align}
\begin{split}
\langle R_s, f \rangle &= R_s(0)f(0)+R'_s(0)f'(0)+\sum_{i=0}^{n}w_i\int_{t_i}^{t_{i+1}}R''_s(u)f''(u)du\\
&=f(0)+sf'(0)+\sum_{i=0}^{n}w_i \int_{t_i}^{t_{i+1}}  \sum_{j=0}^n w_j^{-1} \int_{t_j}^{t_{j+1}} (s-u)_+\delta(v-u)du f''(v)dv\\
&= f(0)+sf'(0)+\sum_{i=0}^n\int_{t_i}^{t_{i+1}}(s-u)_+f''(u)du \\
&=f(s)
\end{split} \\
\begin{split}
\langle \dot{R}_s, f \rangle &= \dot{R}_s(0)f(0)+\dot{R}'_s(0)f'(0)+\sum_{i=0}^{n}w_i\int_{t_i}^{t_{i+1}}\dot{R}_s''(u)f''(u)du\\
&=f'(0)+\sum_{i=0}^{n}w_i \int_{t_i}^{t_{i+1}}  \sum_{j=0}^n w_j^{-1} \int_{t_j}^{t_{j+1}} \Theta(s-u)\delta(v-u)du f''(v)dv\\
&=f'(0)+\sum_{i=0}^n\int_{t_i}^{t_{i+1}}\Theta(s-u)f''(u)du \\
&=f'(s)
\end{split}
\end{align}

Define the two terms of the reproducing kernel $R(s,t)=R_s(t)=R_0(s,t)+R_1(s,t)$, where
\begin{align}
R_0(s,t)&=1+st \\
R_1(s,t)&=\sum_{i=0}^{n}w_i^{-1}\int_{t_i}^{t_{i+1}}(s-u)_+(t-u)_+du
\end{align}
are both non-negative definite themselves. For $R_0$ there corresponds the space of polynomials $\mathcal{H}_0=\lbrace f:f''=0\rbrace$ with an inner product $\langle f,g\rangle = f(0)g(0)+f'(0)g'(0)$, and for $R_1$ there corresponds a sequence of orthogonal spaces $\mathcal{H}^{(i)}$
\begin{equation*}
\mathcal{H}^{(i)} = \lbrace f:f(0)=0,f'(0)=0,\int_{t_i}^{t_{i+1}}f''(t)^2dt<\infty \rbrace 
\end{equation*}
and $\mathcal{H}_1=\oplus _{i=1}^{n-1}\mathcal{H}^{(i)}$. The inner product through the entire space $\mathcal{H}_1$ is $\langle f,g \rangle = \sum_{i=1}^{n-1}w_i\int_{t_i}^{t_{i+1}}f''(t)g''(t)dt$.

Given a sequence of paired sampling points $\lbrace s_i,y_i,v_i \rbrace, i=1,\ldots,n$ on the interval $\lbrack s_1,s_n\rbrack$, it can be transformed to $\lbrace t_i,y_i,v_i \rbrace$ on the interval $\lbrack 0,1\rbrack$, where $0=t_0<t_1<\cdots <t_n<t_{n+1}=1$. The objective function of a V-spline on $\lbrack 0,1\rbrack$ is  
\begin{equation}\label{maineq2}
J[f]=\frac{1}{n}\sum_{i=1}^{n}\left(y_i-f(t_i)\right)^2+\frac{\gamma}{n}\sum_{i=1}^{n}\left(v_i-f'(t_i)\right)^2+\sum_{i=0}^{n}\lambda_i\int_{t_i}^{t_{i+1}}f''(t)^2dt.
\end{equation}
Any $f\in\mathcal{C}_{\mbox{\scriptsize p.w.}}^{(2)}[0,1]$ can be written as 
\begin{equation}\label{GaussianProcessFunctionF2}
f(t)=d_1+d_2t+\sum_{j=1}^{n}c_jR_1(t_j,t)+\sum_{j=1}^{n}b_j\dot{R}_1(t_j,t) +\rho(t)
\end{equation}
Thus, by substituting to the equation \eqref{maineq2}, it can be written as 
\begin{align}\label{GassianProcessRawequation2}
\begin{split}
nJ[f]=& \left( y_i - d_1-d_2t_i-\sum_{j=1}^{n}c_jR_1(t_j,t_i)-\sum_{j=1}^{n}b_j\dot{R}_1(t_j,t_i)-\rho(t_i) \right) ^2\\
+&\gamma\sum_{i=1}^{n} \left( v_i - d_2-\sum_{j=1}^{n}c_jR'_1(t_j,t_i)-\sum_{j=1}^{n}b_j\dot{R}'_1(t_j,t_i)-\rho'(t_i) \right) ^2\\
+&n \sum_{i=0}^{n}\lambda_i \int_{t_i}^{t_{i+1}}  \left( \sum_{j=1}^{n}c_jR''_1(t_j,t)+\sum_{j=1}^{n}b_j\dot{R}''_1(t_j,t)+\rho''(t)\right)^2dt. 
\end{split}
\end{align}
Because of orthogonality, $\rho(t_i) = \langle R_1(t_i,\cdot),\rho\rangle=0$, $\rho'(t_i) = \langle \dot{R}_1(t_i,\cdot),\rho'\rangle=0$, $i=1,\ldots,n$. 
For further use, we need to notice the property of the inner product and $R_1$ satisfy 
\begin{align}
\langle R_1(s,\cdot),\dot{R}_1(t,\cdot)\rangle = R'_1(s,t) \\
\langle \dot{R}_1(s,\cdot),\dot{R}_1(t,\cdot)\rangle = \dot{R}'_1(s,t) 
\end{align}

By denoting the matrices $\{S\}_{ij}=(t_i)^{j-1}$, $j=1,2$, $\{Q\}_{ij}=R_1(t_j,t_i)$, $\{P\}_{ij}=\dot{R}_1(t_j,t_i)$ and $\{P'\}_{ij}=\dot{R}'_1(t_j,t_i)$, the above equation \eqref{GassianProcessRawequation2} becomes the matrix form 
\begin{equation}\label{matriteq2}
\begin{split}
nJ[f]=&\left(\mathbf{y}-S\mathbf{d}-Q\mathbf{c}-P\mathbf{b}\right)^\top \left(\mathbf{y}-S\mathbf{d}-Q\mathbf{c}-P\mathbf{b}\right)\\
+&\gamma\left(\mathbf{v}-S'\mathbf{d}-Q'\mathbf{c}-P'\mathbf{b}\right)^\top \left(\mathbf{v}-S'\mathbf{d}-Q'\mathbf{c}-P'\mathbf{b}\right)\\
+&n\Lambda\left(\mathbf{c}^\top Q\mathbf{c} + 2\mathbf{c}^\top P\mathbf{b}+ \mathbf{b}^\top P'\mathbf{b}\right)+n\Lambda(\rho,\rho),
\end{split}
\end{equation}
where $\lambda_i=\Lambda w_i$. 

Thus, the solution to \eqref{maineq2} is computed via the minimization of the first three terms in \eqref{matriteq2} with respect to $\mathbf{d}$, $\mathbf{c}$ and $\mathbf{b}$.

Therefore, the calculation goes through the same process in Section \ref{sectionBayesEstimate} and the following theorem is obtained. 
\begin{theorem}
The smoothing V-spline of \eqref{maineq2} is the posterior mean of $f=f_0+f_1 + \dot{f}_1$, where $f_0$ diffuses in span $\left\lbrace 1,t\right\rbrace$ and $f_1$, $\dot{f}_1$ have Gaussian process priors with mean zero and covariance functions
\begin{align}
\Cov \left(f_1,f_1\right)   &= \beta R_1\left(s,t\right)  =\beta \sum_{i=0}^n w_i^{-1} \int_{t_i}^{t_{i+1}} \left(s-u\right)_+\left(t-u\right)_+du, \\
\Cov \left(\dot{f}_1,f_1\right)  &= \beta \dot{R}_1\left(s,t\right) =\beta \sum_{i=0}^n w_i^{-1} \int_{t_i}^{t_{i+1}} \Theta\left(s-u\right)\left(t-u\right)_+du,\\
\Cov \left(\dot{f}_1,\dot{f}_1\right)  &= \beta \dot{R}'_1\left(s,t\right) =\beta \sum_{i=0}^n w_i^{-1} \int_{t_i}^{t_{i+1}} \Theta\left(s-u\right)\Theta\left(t-u\right)du,
\end{align}
for $\beta = \sigma^2/n\Lambda$.
\end{theorem}

\section{V-Spline with Correlated Random Errors}

In most of the studies on polynomial smoothing splines, the random errors are assumed being independent. By contrast,  observations are often correlated in applications, such as time series data and spatial data. It is known that the correlation greatly affects the selection of smoothing parameters, which are critical to the performance of smoothing spline estimates \citep{wang1998smoothing}.  The parameter selection methods, such as generalized maximum likelihood (GML), generalized cross-validation (GCV), underestimate smoothing parameters when data are correlated.

\cite{diggle1989spline} extend GCV for choosing the degree of smoothing spline to accommodate an autocorrelated error sequence, by which the smoothing parameter and autocorrelation parameters are estimated simultaneously.  \cite{kohn1992nonparametric} propose an algorithm to evaluate the cross-validation functions, whose autocorrelated errors are modeled by an autoregressive moving average. \cite{wang1998smoothing} extend GML and unbiased risk (UBR), other than GCV, to estimate the smoothing parameters and correlation parameters simultaneously. In this section, we explore the extended GCV for V-spline with correlated errors.

First of all, consider observations $y=f(t)+\varepsilon_1$ and $v=f'(t)+\varepsilon_2$, where $\varepsilon_1\sim N\left(0,\sigma^2W^{-1}\right)$, $\varepsilon_2\sim N\left(0,\frac{\sigma^2}{\gamma}U^{-1}\right)$ with variance parameter $\sigma^2$, and the structures of correlation matrices $W$ and $U$ are known. The V-spline $\hat{f}$ with correlated errors in the space $\mathcal{C}_{\mbox{\scriptsize p.w.}}^{(2)}[0,1]$ is the minimizer of 
\begin{equation}
\frac{1}{n}\left(\mathbf{y}-\mathbf{f}\right)^\top W\left(\mathbf{y}-\mathbf{f}\right)+\frac{\gamma}{n}\left(\mathbf{v}-\mathbf{f}'\right)^\top U\left(\mathbf{v}-\mathbf{f}'\right)+\lambda\int_0^1\left(f''\right)^2dt.
\end{equation}
Because of $f=\sum_{i=1}^{2n}\theta_iN_i\left(t\right)$ is a linear combination of basis functions $\left\lbrace N_i(t)\right\rbrace_{i=1}^{2n}$, extended to the solution with covariance matrices, the coefficients is found by 
\begin{equation}
\hat{\theta}=\left(B^\top W B+ \gamma C^\top UC+n\Omega_\lambda\right)^{-1}\left(B^\top W \mathbf{y}+\gamma C^\top U\mathbf{v}\right).
\end{equation}
Furthermore, in Gaussian process regression, the covariance matrix with correlated variances becomes \\
$M=\begin{bmatrix}
Q+n\lambda W& P\\
Q'& P'+\frac{n\lambda}{\gamma}U
\end{bmatrix}$ 
and the rest stays the same.

Additionally, it is known that the parameter $\hat{\theta}=\left(B^\top B+\gamma C^\top C+n\Omega_\lambda\right)^{-1}\left(B^\top\mathbf{y}+\gamma C^\top\mathbf{v}\right)$ and will give us the following form \small
\begin{equation}
\begin{split}
 \hat{\mathbf{f}}&=B\hat{\theta}=B\left(B^\top B+\gamma C^\top C+n\Omega_\lambda\right)^{-1}B^\top\mathbf{y}+B\left(B^\top B+\gamma C^\top C+n\Omega_\lambda\right)^{-1} C^\top\mathbf{v}\\&=S\mathbf{y}+\gamma T\mathbf{v},
 \end{split}
 \end{equation}
 \begin{equation}
 \begin{split}
\hat{\mathbf{f}}'&=C\hat{\theta}=C\left(B^\top B+\gamma C^\top C+n\Omega_\lambda\right)^{-1}B^\top\mathbf{y}+C\left(B^\top B+\gamma C^\top C+n\Omega_\lambda\right)^{-1}C^\top \mathbf{v}\\&=U\mathbf{y}+\gamma V\mathbf{v}.
 \end{split}
\end{equation}\normalsize

\begin{lemma}\label{tractorsplinecvscore}
The cross-validation score of a V-spline satisfies
\begin{equation}\label{tractorcv}
\mbox{CV}\left(\lambda,\gamma\right)=\frac{1}{n}\sum_{i=1}^{n} \left( \frac{\hat{f}(t_i)-y_i+\gamma \frac{T_{ii}}{1-\gamma V_{ii}}(\hat{f}'(t_i)-v_i)}{1-S_{ii}-\gamma\frac{T_{ii}}{1-\gamma V_{ii}}U_{ii}} \right)^2
\end{equation}
where $\hat{f}$ is the V-spline smoother calculated from the full data set $\left\lbrace (t_i,y_i,v_i)\right\rbrace$ with smoothing parameter $\lambda$ and $\gamma$.
\end{lemma}

Followed by the approximation $S_{ii}\approx\frac{1}{n}\tr(S)$, $T_{ii}\approx\frac{1}{n}\tr(T)$, $U_{ii}\approx\frac{1}{n}\tr(U)$ and $V_{ii}\approx\frac{1}{n}\tr(V)$ \citep{syed2011review}, the GCV for the V-spline will be 
\begin{equation}
\mbox{GCV}(\lambda,\gamma)=\frac{1}{n}\sum_{i=1}^{n}\left( \frac{\hat{f}(t_i)-y_i+ \frac{\gamma\tr(T)/n}{1-\gamma \tr(V)/n}(\hat{f}'(t_i)-v_i)}{1-\tr(S)/n-\frac{\gamma\tr(T)/n}{1-\gamma \tr(V)/n} \tr(U)/n} \right)^2,
\end{equation}
which may provide further computational savings since it requires finding the trace rather than the individual diagonal entries of the hat matrix. Hence, it can be written in the form of \footnotesize
\begin{equation}
\mbox{GCV}(\lambda,\gamma)=\frac{\left(\mathbf{\hat{f}}-\mathbf{y}\right)^\top \left(\mathbf{\hat{f}}-\mathbf{y}\right) + \frac{2\tr\left(\gamma T\right)}{\tr\left(I-\gamma V\right)}\left(\mathbf{\hat{f}}-\mathbf{y}\right)^\top \left(\mathbf{\hat{f}}'-\mathbf{v}\right) + \left( \frac{\tr(\gamma T)}{\tr(I-\gamma V)} \right)^2 \left(\mathbf{\hat{f}}'-\mathbf{v}\right)^\top \left(\mathbf{\hat{f}}'-\mathbf{v}\right)}{\left( \tr(I-S-\frac{\tr(\gamma T)}{\tr(I-\gamma V)}U) \right)^2}.
\end{equation} \normalsize
A natural extension to the above GCV for V-spline with correlated errors is \scriptsize
\begin{equation}
\mbox{GCV}\left(\lambda,\gamma\right)=\frac{\left(\mathbf{\hat{f}}-\mathbf{y}\right)^\top W\left(\mathbf{\hat{f}}-\mathbf{y}\right)+\frac{2\tr\left(\gamma T\right)}{\tr\left(I-\gamma V\right)}\left(\mathbf{\hat{f}}-\mathbf{y}\right)^\top W^{1/2}U^{\top 1/2}\left(\mathbf{\hat{f}}'-\mathbf{v}\right) + \left( \frac{\tr\left(\gamma T\right)}{\tr\left(I-\gamma V\right)} \right)^2 \left(\mathbf{\hat{f}}'-\mathbf{v}\right)^\top U\left(\mathbf{\hat{f}}'-\mathbf{v}\right)}{\left( \tr\left(I-S-\frac{\tr\left(\gamma T\right)}{\tr\left(I-\gamma V\right)}U\right) \right)^2}.
\end{equation}
\normalsize

The GCV is used for finding the unknown constant parameter $\lambda$, instead of a piecewise constant $\lambda(t)$ at different intervals, and the parameter $\gamma$. If the errors are independent, in which way $W$ and $U$ become identity matrices, the solution $\hat{f}$ degenerates to a conventional V-spline with constant $\lambda$ through over the entire interval $[0,1]$.  

For a non-trivial  V-spline, the parameter $\lambda$ can be easily substituted by $\Lambda$, which then is optimized by the above formula.

\section{Conclusion}

In this paper, we discussed the correspondence between polynomial smoothing spline and Bayes estimates given improper priors. In fact, the smoothing spline is a particular case of Gaussian process regression. By following the work done by \cite{gu2013smoothing}, we find the Bayes estimate for V-splines in two scenarios: constant penalty parameters $\lambda$ and $\gamma$ on the entire interval $\lbrack0,1\rbrack$; parameters $\Lambda$ and $\gamma$, but $\Lambda$ is distributing on different subintervals and its value depending on the weight of that subinterval. Additionally, we give the formula of GCV for V-spline with correlated errors on $y$ and $v$.

\bibliographystyle{Chicago}
\bibliography{/Users/zcao/Documents/Latex/Working/WorkingAll}

\end{document}